\documentclass[12pt]{article}

\usepackage{amsmath}
\usepackage{amssymb}
\usepackage{url}

\newtheorem{Thm}{Theorem}
\newtheorem{Lem}[Thm]{Lemma}
\newenvironment{proof}{\noindent \emph{Proof.}}{\hspace*{\fill} $\Box$\newline}
\newenvironment{Proof}{\begin{proof}}{\end{proof}}

\newcommand{\LF}{\ensuremath{\mathrm{LF}}}
\newcommand{\NF}{\ensuremath{\mathrm{NF}}}
\newcommand{\Aut}{\ensuremath{\mathrm{Aut}}}

\newcommand{\onef}[1]{\ensuremath{\mathcal{F}({#1})}}
\newcommand{\ext}{\ensuremath{\mathrm{ext}}}
\newcommand{\seeds}{\ensuremath{\mathrm{seeds}}}

\begin{document}

%%%%%%%%%%%%%%%%%%%%%%%%%%%%%%%%%%%%%%%%%%%%%%%%%% Front matter information %%%

\title{There are 1,132,835,421,602,062,347\\ nonisomorphic 
one-factorizations of $K_{14}$}
\author{Petteri Kaski\\
Helsinki Institute for Information Technology HIIT\\
University of Helsinki, Department of Computer Science\\
P.O.\ Box 68, 00014 University of Helsinki, Finland\\
E-mail: {\tt petteri.kaski@cs.helsinki.fi}
\and
Patric R. J. \"Osterg\aa rd\\
Department of Communications and Networking\\
Helsinki University of Technology\\
P.O.\ Box 3000, 02015 TKK, Finland\\
E-mail: {\tt patric.ostergard@tkk.fi}}
\date{}

\maketitle

%%%%%%%%%%%%%%%%%%%%%%%%%%%%%%%%%%%%%%%%%%%%%%%%%%%%% Abstract and keywords %%%

\begin{abstract}
We establish by means of a computer search 
that a complete graph on 14 vertices has 
98,758,655,816,833,727,741,338,583,040 distinct and
1,132,835,421,602,062,347 nonisomorphic
one-factorizations. The enumeration is constructive for 
the 10,305,262,573 isomorphism classes that admit a nontrivial 
automorphism.
\end{abstract}

%%%%%%%%%%%%%%%%%%%%%%%%%%%%%%%%%%%%%%%%%%%%%%%%%%%%%%%%%%%%% Document body %%%

\section{Introduction}

A \emph{one-factor} of a graph is a 1-regular spanning subgraph.
A \emph{one-factorization} of a graph is a set of one-factors
such that every edge of the graph occurs in a unique one-factor.
Two one-factorizations are {\em isomorphic} if there is a bijection 
between the vertex sets of the underlying graphs that maps 
the one-factors in one one-factorization onto the one-factors of 
the other. One-factorizations of complete graphs are basic
combinatorial objects with a wide variety of applications \cite{A,W1,W2}.

A fundamental problem for any family of combinatorial objects is that
of classifying these for small parameters, that is, of listing exactly
one representative from each isomorphism class. The complete graph $K_n$ on $n$ vertices has, up to isomorphism,
a unique one-factorization for $n=2,4,6$. The graph $K_8$ has 6 
nonisomorphic one-factorizations; this result was 
obtained by Dickson and Safford \cite{DS} almost exactly one century ago. 
In 1974, Gelling and Odeh \cite{GO} published the result that $K_{10}$ has
396 nonisomorphic one-factorizations. 
A result by Dinitz, Garnick,
and McKay \cite{DGM} showing that there are 526,915,620 nonisomorphic 
one-factorizations of $K_{12}$ appeared exactly two decades after
the result by Gelling and Odeh. The result for $K_{12}$ was one of the most 
extensive computer classifications at that time, and---with an estimated 
$1.132\times 10^{18}$ nonisomorphic one-factorizations \cite{DGM} of 
$K_{14}$---the current limit for methods that rely on constructively listing 
one specimen from each isomorphism class was thereby reached.

For many types of combinatorial objects, no essentially faster 
method than explicitly constructing one representative from each 
isomorphism class is known for counting the number of isomorphism 
classes. However, in cases where one is able to count the ``labeled'' 
objects substantially faster, the following alternative approach
for determining the number of isomorphism classes becomes possible.

Let $\Gamma$ be a finite group that acts on a finite set $\Omega$.
Denote by $N_i$ the number of orbits on $\Omega$ whose elements
have stabilizer subgroups of order $i$ in $\Gamma$. Then, by the
Orbit-Stabilizer Theorem,
\begin{equation}
\label{eq:totaln}
|\Omega| = |\Gamma|\sum_i \frac{N_i}{i}.
\end{equation}
If we know $|\Gamma|$, $|\Omega|$, and $N_i$ for each $i \geq 2$, then
we can solve \eqref{eq:totaln} for $N_1$ and 
thereby obtain the number of orbits $\sum_i N_i$.
This observation can be applied in a combinatorial context
by considering a group action for which 
(i) the set $\Omega$ consists of the ``labeled'' objects;
(ii) the orbits of $\Gamma$ on $\Omega$ correspond to the isomorphism 
     classes of objects (the ``unlabeled'' objects); and 
(iii) the stabilizer subgroup of an object under the action of $\Gamma$ 
      corresponds to the automorphism group of the object. 
In particular, the values $N_i$ for $i\geq 2$ are obtained by classifying 
up to isomorphism the objects that have a nontrivial automorphism group.

Latin squares and one-factorizations of complete graphs are two types
of objects for which the outlined idea is effective. Recently, 
McKay, Meynert, and Myrvold \cite{MMM} successfully applied this general 
idea to the problem of enumerating Latin squares. In the present paper, 
it is applied to the problem of enumerating one-factorizations of 
the complete graph $K_{14}$. 

The approach for counting the number of distinct one-factorizations of 
$K_{14}$ is discussed in Section~\ref{sect:distinct}; there are 
98,758,655,816,833,727,741,338,583,040 such objects. Classification
of one-factorizations of $K_{14}$ with a nontrivial automorphism group is 
considered in Section~\ref{sect:nonisomorphic}. 
Application of \eqref{eq:totaln}
then reveals that there are 1,132,835,421,602,062,347
isomorphism classes of one-factorizations of $K_{14}$. 
Particular emphasis is put on verifying the correctness of these 
computational results.

{\em Conventions.}
For standard graph-theoretic terminology
we refer to \cite{W3}. All graphs considered in this paper 
are undirected and without multiple edges or loops. 
For a graph $G$, denote by $\bar G$ the complement of $G$,
by $\onef G$ the set of all one-factors of $G$, 
by $[G]$ the isomorphism class of $G$, by $\Aut(G)$ the automorphism
group of $G$, by $\LF(G)$ the number of distinct one-factorizations 
of $G$, and by $\NF(G)$ the number of nonisomorphic one-factorizations of $G$.

\section{Counting distinct one-factorizations}
\label{sect:distinct}

To count the number of distinct one-factorizations of $K_{14}$, 
that is, to compute $\LF(K_{14})$, we essentially rely on a recursion 
pioneered by Dinitz, Garnick, and McKay \cite[Sect.~4]{DGM}, who
used the recursion to check their constructive enumeration of 
the isomorphism classes of one-factorizations of $K_{12}$. 

In our case, however, we do not have an independently computed 
value of $\LF(K_{14})$ available. 
To gain confidence that the computed value of $\LF(K_{14})$ 
is correct, we will use 
(i) a new ``forward accumulation'' technique to compute 
    the value of $\LF(K_{14})$; and 
(ii) Dinitz--Garnick--McKay recursion to \emph{check} 
    both the intermediate results and the final computed value.

It is convenient to start by describing Dinitz--Garnick--McKay 
recursion and then proceed to describe the new technique.

\subsection{Dinitz--Garnick--McKay recursion}

Consider a graph $G$ that has a one-factorization.
Such a graph is necessarily $k$-regular for some nonnegative
integer $k$. Avoiding trivial cases, we assume that $k\geq 2$.
Let $F$ be a one-factor in $G$, and denote by $G-F$ 
the $(k-1)$-regular graph obtained by deleting the edges 
of $F$ from $G$. To arrive at a recursion, suppose that we know 
the value of $\LF(G-F)$ for every one-factor $F$ of $G$. 
By counting in two different ways the number of distinct 
one-factorizations of $G$ with one individualized one-factor, 
we have
\begin{equation}
\label{eq:dgm}
k\cdot\LF(G)=\sum_{F\in\onef G} \LF(G-F).
\end{equation}
In particular, if we know $\LF(H)$ for every $(k-1)$-regular graph $H$, then
we can compute $\LF(G)$ for any given $k$-regular graph $G$ via \eqref{eq:dgm}.
This is \emph{Dinitz--Garnick--McKay recursion} from
$(k-1)$-regular to $k$-regular graphs.

In practice it suffices to evaluate the right-hand side of 
\eqref{eq:dgm} only for one-factors that contain a fixed edge 
of $G$, in which case the multiplication by $k$ on the left-hand
side is not necessary. Similarly, $\LF(H)$ needs to be computed
for only one graph $H$ in each isomorphism class of regular graphs. 

Dinitz--Garnick--McKay recursion has the property that
it considers each $k$-regular graph in turn, and ``looks back''
at the $(k-1)$-regular graphs. 
We describe next a new technique that ``looks forward'' 
at the $k$-regular graphs from each $(k-1)$-regular graph in turn.

\subsection{Forward accumulation}

The forward accumulation approach is based on the
immediate observation that every $k$-regular 
one-factorizable graph $G$ can be decomposed into 
a union $G=H\cup F$ of 
(i) a $(k-1)$-regular one-factorizable graph $H$; and 
(ii) a one-factor $F\in\onef{\bar H}$.
Put otherwise, we can visit every isomorphism class of 
$k$-regular one-factorizable graphs by the following
procedure: for each isomorphism class $[H]$ 
of $(k-1)$-regular one-factorizable graphs, consider exactly one graph $H$ 
from the isomorphism class; for each such graph $H$, consider each
one-factor $F\in\onef{\bar H}$; for each such pair $(H,F)$, 
visit the isomorphism class $[H\cup F]$. 

To compute the value $\LF(G)$ for each visited isomorphism
class $[G]$, we associate with $[G]$ an accumulator
variable $x_{[G]}$ that is initially set to zero
and incremented whenever $[G]$ is visited. 
Our objective is to have the value $k\cdot\LF(G)$ in
the accumulator when the visiting procedure halts.
To determine an appropriate increment to $x_{[G]}$ on each visit, 
we proceed to analyze the visiting procedure in more detail.

To this end, consider the set of \emph{all} pairs $(H,F)$ 
such that $H$ is a $(k-1)$-regular one-factorizable graph 
and $F\in\onef{\bar H}$. Let us view two such pairs as
\emph{isomorphic} if one can be obtained from the other by 
relabeling the vertices. 
The following lemmata are immediate consequences of 
the Orbit-Stabilizer Theorem.

\begin{Lem}
\label{lem:decomp}
Any graph $G$ in the class $[H\cup F]$ admits exactly 
\[
\sigma(H,F)=\frac{|\Aut(H\cup F)|}{|\Aut(H)\cap\Aut(F)|}
\]
decompositions $G=H'\cup F'$ 
into pairs $(H',F')$ in the class $[(H,F)]$.
\end{Lem}

\begin{Lem}
\label{lem:visits}
The procedure visits a class $[H\cup F]$ exactly 
\[
\tau(H,F)=\frac{|\Aut(H)|}{|\Aut(H)\cap\Aut(F)|}
\]
times via pairs $(H',F')$ in the class $[(H,F)]$.
\end{Lem}

It now follows from Lemma \ref{lem:decomp} 
and \eqref{eq:dgm} that, for any $k$-regular graph $G$,
\begin{equation}
\label{eq:dgm-up}
\sum_{[(H,F)]:H\cup F=G}\sigma(H,F)\cdot\LF(H)=
\sum_{(H,F):H\cup F=G}\LF(H)=
k\cdot \LF(G).
\end{equation}
This observation enables us to accumulate the value $k\cdot\LF(G)$ 
for each $k$-regular $G$.
Namely, each time $[G]$ is visited via a pair $(H,F)$, 
we increment $x_{[G]}$ by the rule
\begin{equation}
\label{eq:acc}
x_{[G]}\leftarrow x_{[G]}+\sigma(H,F)\cdot\tau(H,F)^{-1}\cdot\LF(H).
\end{equation}
Equivalently, for each pair $(H,F)$ considered by the visiting procedure,
we apply the rule
\begin{equation}
\label{eq:acc2}
x_{[H\cup F]}\leftarrow x_{[H\cup F]}+
                        \frac{|\Aut(H\cup F)|}{|\Aut(H)|}\cdot\LF(H).
\end{equation}

\begin{Lem}
The total accumulation to $x_{[G]}$ is $k\cdot \LF(G)$.
\end{Lem}
\begin{Proof}
By Lemma \ref{lem:visits} and \eqref{eq:acc}, 
the total accumulation to $x_{[G]}$ from a class $[(H,F)]$ satisfying 
$[H\cup F]=[G]$ is $\sigma(H,F)\cdot\LF(H)$. Taking the sum over all
such classes, the claim follows by \eqref{eq:dgm-up}.
\end{Proof}

\subsection{Implementation details}

Starting with the empty graph of order 14, that is, the
unique 0-regular graph on 14 vertices, we use the forward accumulation
technique for each $k=1,2,\ldots,13$ in turn to compute both
(i) exactly one graph $G$ from each isomorphism class 
    of $k$-regular one-factorizable graphs; and
(ii) the value $\LF(G)$ for the graphs $G$.

The representative graph $G$ in an isomorphism class is the canonical
form computed by \emph{nauty} (version 2.2) \cite{M1} using the
built-in \texttt{adjtriang} vertex invariant. The canonical form 
is stored in 16 bytes of memory as a bit map of the upper triangle of 
the adjacency matrix; the associated accumulator variable uses 
32 bytes. To enable rapid searching of these 48-byte records, 
we use an open-addressing hash table with $10^9$ entries 
indexed by 4-byte hash values of the 16-byte bit maps. 

We use the GNU Multiple Precision Arithmetic Library \cite{gmplib}
to carry out the accumulator arithmetic.
An elementary backtrack algorithm suffices for listing the
one-factors $F\in\onef{\bar H}$.

The number of isomorphism classes of $k$-regular one-factorizable 
graphs on 14 vertices is, for $k=0,1,\ldots,13$,
\[
1, 
1, 
4, 
504, 
87977, 
3459360, 
21609293,
21609301, 
3459386, 
88193, 
540, 
13, 
1, 
1.
\]
The performance bottleneck is at $k=7$, where 43,218,594
records need to be stored, occupying about
6 GB of memory together with the hash table (cf.\ \cite{MW}). 
In terms of running time, the entire computation 
(including the correctness checks described in what follows) 
took about 13 days on a Linux PC with a 3.66-GHz Intel Xeon CPU
and 32 GB of main memory.

\subsection{Correctness checks}

Based on forward accumulation, we have the value
$\LF(G)$ for every regular graph $G$ on 14 vertices.
In particular, 
\begin{equation}
\label{eq:lfk14}
\LF(K_{14})=\text{98,758,655,816,833,727,741,338,583,040}.
\end{equation}

As a first check, we use Dinitz--Garnick--McKay 
recursion to verify the $\LF(G)$ values for each regular $G$.
The same values are obtained both using forward accumulation
and Dinitz--Garnick--McKay recursion.

As a second check, we use Meringer's classification program 
{\bf genreg} \cite{Me} to generate all the regular graphs on 14
vertices, and then filter out the graphs that do not have 
a one-factorization. The obtained graphs agree with those
obtained by forward accumulation.

As a third check, we use the following observation due to 
Dinitz, Garnick, and McKay.
Taking the sum over all isomorphism classes $[G]$ of $k$-regular
graphs on $n$ vertices, we have
\begin{equation}
\label{eq:meet-in-the-middle}
\LF(K_n)=
\binom{n-1}{k}^{-1}\sum_{[G]}\frac{n!}{|\Aut(G)|}\cdot\LF(G)\cdot\LF(\bar G).
\end{equation}
Moreover, this holds for every $k=0,1,2,\ldots,n-1$.
Using data obtained from forward accumulation, we evaluate
the right-hand side of \eqref{eq:meet-in-the-middle} for $n=14$
and $k=0,1,\ldots,13$. In each case the computed value agrees
with \eqref{eq:lfk14}.

To enable further checks, we display in Table \ref{tbl:lf-k-reg}
the value 
\[
\sum_{[G]}\frac{14!}{|\Aut(G)|}\cdot\LF(G)
\]
for each $k=0,1,\ldots,13$, where the sum is taken over all 
isomorphism classes $[G]$ of $k$-regular graphs on 14 vertices. 
Put otherwise, the tabulated value
is the number of distinct one-factorizations of $k$-regular graphs
on a fixed set of 14 vertices; cf.~\cite[Table~7]{DGM}, where a slightly
different quantity is tabulated, however.

\begin{table}
\begin{center}
\begin{tabular}{rr}\hline
\raisebox{-0.5mm}{$k$} & \raisebox{-0.5mm}{Distinct one-factorizations}\\[0.5mm]\hline
\raisebox{-0.5mm}{1}  & \raisebox{-0.5mm}{135135}\\[0.5mm]
2  &                      5338373040\\
3  &                  78634135419840\\
4  &              461142306338313600\\
5  &          1078882420304271623040\\
6  &        972197327694773750169600\\
7  &     315828427387711768964628480\\
8  &   33491835583595013396417085440\\
9  & 1006698095378044123991615078400\\
10 & 7024525682952576878777802424320\\
11 & 8573318527281503086919968358400\\
12 & 1283862525618838460637401579520\\
13 &   98758655816833727741338583040\\
\hline
\end{tabular}
\end{center}
\bigskip
\caption{Distinct one-factorizations of $k$-regular graphs on 14 vertices}
\label{tbl:lf-k-reg}
\end{table}

\section{Counting nonisomorphic one-factorizations}
\label{sect:nonisomorphic}

Our primary objective in this section is to classify
the one-factorizations of $K_{14}$ with nontrivial automorphisms.
Once this classification is available, it is a simple matter to
determine $\NF(K_{14})$ based on $\LF(K_{14})$ and the classification
data. Before presenting the classification approach, we discuss
a representation for one-factorizations of the complete graph in the
framework of group divisible designs, and narrow down the
automorphisms that need to be considered to obtain a complete
classification.

\subsection{One-factorizations as group divisible designs}
\label{sect:gdd}

It is well known that a one-factorization of $K_n$ can be viewed 
as a particular type of group divisible design (GDD). 
For our purposes this representation will be 
convenient as it enables us to immediately prescribe the action of 
an automorphism not only on the vertices but also on 
the one-factors.

To develop the GDD representation, 
let $U=\{u_1,u_2,\ldots,u_{n-1}\}$ be a set with 
one element for each one-factor, and let $V=\{v_1,v_2,\ldots,v_n\}$ be 
a set (disjoint from $U$) with one element for each vertex of 
$K_n$. The elements of $U\cup V$ are called \emph{points}.
We can now let a 3-subset of the form $\{u_k,v_i,v_j\}$
carry the information that the edge $\{v_i,v_j\}$ occurs in the 
one-factor $u_k$. In this setting, a \emph{one-factorization} 
of $K_n$ is a tuple $\mathcal{X}=(U,V,\mathcal{B})$, where $\mathcal{B}$ 
is a set of 3-subsets of points, called \emph{blocks}, such that

\begin{itemize}
\item[(a)] for all $1\leq k\leq n-1$ and $1\leq i\leq n$,
the pair $\{u_k,v_i\}$ occurs in a unique block; 
\item[(b)] for all $1\leq i<j\leq n$, 
the pair $\{v_i,v_j\}$ occurs in a unique block; and
\item[(c)] for all $1\leq k<\ell\leq n-1$, no block contains the
pair $\{u_k,u_\ell\}$. 
\end{itemize}

A reader acquainted with GDDs immediately observes that 
$\mathcal{X}$ is in fact a GDD with group type $(n-1)^11^n$, 
block size $k=3$, and index $\lambda=1$; however, here we find 
it more convenient to speak of one-factorizations.

Two one-factorizations, $\mathcal{X}=(U,V,\mathcal{B})$ and 
$\mathcal{X}'=(U',V',\mathcal{B}')$, are \emph{isomorphic} if there exists 
a bijection $\varphi:U\cup V\rightarrow U'\cup V'$ such that
$\varphi(U)=U'$, $\varphi(V)=V'$, and $\varphi(\mathcal{B})=\mathcal{B}'$.
Such a $\varphi$ is an \emph{isomorphism} from $\mathcal{X}$ onto 
$\mathcal{X}'$. An isomorphism of $\mathcal{X}$ onto itself is an 
\emph{automorphism} of $\mathcal{X}$. We denote by $\Aut(\mathcal{X})$
the automorphism group of $\mathcal{X}$.
Note that the restriction of an isomorphism 
$\varphi$ to $V$ uniquely determines $\varphi$ on $U$. 
It follows that the standard (graphical) and the
GDD representations of one-factorizations of $K_n$ are 
equivalent for purposes of classification up to isomorphism.

\subsection{Automorphisms of one-factorizations}
\label{sect:autotype}

Any nontrivial group has a subgroup of prime order, and our
classification considers all possible
such groups that can occur as a group of automorphisms 
of a one-factorization.
Ihrig and Petrie \cite{IP} carried out an extensive study of
\emph{all} possible automorphisms for one-factorizations---in
particular, for the complete graph $K_{12}$---but
we indeed only need those with prime order.

It is convenient to assume in what follows that the sets $U$ and $V$ 
are arbitrary but fixed. This enables the following two simplifications.
First, isomorphisms between one-factorizations are permutations of 
$U\cup V$ that fix $U$ and $V$ setwise. Denote by $\Gamma$ the group of 
all such permutations of $U\cup V$. Second, we can identify
a one-factorization $\mathcal{X}$ with its set of blocks $\mathcal{B}$.

A group element $\alpha\in\Gamma$ of prime order $p$ is determined 
up to conjugation in $\Gamma$ by the number of fixed points
it has in the sets $U$ and $V$. Denote by $f_U$ and $f_V$ the number
of fixed points of $\alpha$ in $U$ and $V$, respectively.
Because the cycle decomposition of $\alpha$ consists only of 
fixed points and $p$-cycles, it is immediate that $p$ must divide 
both $|U|-f_U=n-1-f_U$ and $|V|-f_V=n-f_V$. 
Not all such types $(p,f_U,f_V)$ define automorphisms of 
one-factorizations, however.

We proceed to narrow down the possible types $(p,f_U,f_V)$. 
The following lemma is analogous to \cite[Lemma 32]{K};
see also \cite[Lemma 4.1]{SS}.

\begin{Lem}
\label{lem:autotype}
Let $\alpha$ be an automorphism of a one-factorization $\mathcal{X}$
with $f_U\geq 1$ and $f_V\geq 1$.
Then, $\mathcal{X}$ restricted to the fixed points of $\alpha$ 
forms a one-factorization. In particular, 
$f_U=f_V-1$ and $f_V$ is even.
\end{Lem}

\begin{Proof}
Let $x$ and $y$ be two distinct points fixed by $\alpha$, at least
one of which is in $V$. Such points clearly exist if $f_U\geq 1$ and 
$f_V\geq 1$. Then there is exactly one block $\{x,y,z\}$ 
of $\mathcal{X}$ that contains both $x$ and $y$. Since $\alpha$ is 
an automorphism of $\mathcal{X}$, also 
$\{x,y,\alpha(z)\}$ is a block, and thus $\alpha(z) = z$.
This shows that no block of $\mathcal{X}$ intersects the set of fixed 
points of $\alpha$ in 2 points. Consequently, every block 
of $\mathcal{X}$ intersects the set of fixed points in 0, 1, or 3 points.
Disregarding all other blocks except those intersecting in 3 points,
we obtain a set system $\mathcal{X}'=(U',V',\mathcal{B}')$,
where $U'$ and $V'$ are the sets of points 
fixed by $\alpha$ in $U$ and $V$, respectively.

To see that $\mathcal{X}'$ is a one-factorization
of a complete graph of order $n'=f_V$, we first count in two different 
ways the tuples $(u_k,v_i,B)$ such that 
$u_k\in U'$, $v_i\in V'$, $B\in\mathcal{B}'$, and $\{u_k,v_i\}\subseteq B$. 
We have $|U'|\cdot |V'| = 2|\mathcal{B}'|$.
As all 2-subsets of $V'$ occur in a unique block, we have
$|\mathcal{B'}| = |V'|(|V'|-1)/2$. 
Combining the two equalities and noticing that $|U'|=f_U$ and
$|V'|=f_V$, we get $f_U = f_V-1$.
Furthermore, it follows that $\mathcal{X}'$ meets the requirements 
(a), (b), and (c) in the definition of a one-factorization of a complete 
graph for $n'=f_V$. In particular, $f_V$ must then be even.
\end{Proof}

The following lemma is due to Seah and Stinson 
\mbox{\cite[Lemmata 4.2 and 4.3]{SS}}.

\begin{Lem}
Any nonidentity automorphism of a one-factorization satisfies 
$f_V \leq n/2$. Equality holds only if the automorphism 
has order $2$.
\end{Lem}

The following two lemmata rule out certain automorphism types with $p=2$; 
the argument in the proof of the second one is similar to that 
of \mbox{\cite[Lemma 33]{K}}.

\begin{Lem}
Let $n \equiv 2 \pmod{4}$.
If $p=2$ and $f_V=0$ for an automorphism of a one-factorization,
then $f_U\leq n/2$.
\end{Lem}

\begin{Proof}
With the objective of bounding $f_U$, consider any one-factor 
$u_k\in U$ fixed by $\alpha$. Since the number of edges 
in the one-factor is $n/2 \equiv 1 \pmod{2}$, 
not all orbits of edges in the one-factor can have size 2. Thus, 
because $f_V=0$, there is at least one block $\{u_k,v_i,v_j\}$ 
such that $(v_i\ v_j)$ is a 2-cycle of $\alpha$. 
The number of such 2-cycles in $\alpha$ is $n/2$, and thereby 
$f_U \leq n/2$.
\end{Proof}

\begin{Lem}
\label{lem:none2}
Let $n \equiv 4 \pmod{8}$ or $n \equiv 6 \pmod{8}$.
If $p=2$ and $f_V=0$ for an automorphism of a one-factorization,
then $f_U\neq 1$.
\end{Lem}

\begin{Proof}
To reach a contradiction, assume that $\alpha$ is an automorphism
of a one-factorization with $p=2$, $f_V=0$, and $f_U=1$.
Without loss of generality we may assume that
$U\cup V=\{-(n-1),-(n-2),\ldots,n-1\}$
and that $\alpha(x)=-x$ for all $x\in U\cup V$.
In particular, because $\alpha$ is an automorphism,
both $U$ and $V$ are closed under negation.
Because no pair of points occurs in more than one block,
it follows by $\alpha(\{x,-y,y\}) = \{\alpha(x),y,-y\}$
that $\alpha(x)=x$ and $x=0$.
In particular, each of the $n/2$ pairs of the form $\{-y,y\}\subseteq V$ 
occurs in a block of the form $\{0,-y,y\}$. 
All the remaining blocks are moved by $\alpha$, that is, $\{x,y,z\}$ 
is a block if and only if $\{-x,-y,-z\}$ is a block, and each such pair 
of blocks contains either 0 or 4 pairs of points with opposite 
signs. Considering pairs of points with either one point in $U$
and one in $V$ or two points in $V$, there are 
$2\cdot(n-2)/2\cdot n/2+(n/2)^2$ such pairs of points
with opposite signs. 
The fixed blocks account for $n/2$ occurrences, so 
the moved blocks must account for $3\cdot n/2\cdot (n/2-1)$ occurrences.
For $n \equiv 4 \pmod{8}$ and $n \equiv 6 \pmod{8}$ this
number is not divisible by $4$, a contradiction.
\end{Proof}

For $n=14$, prime orders $p=2,3,5,7,11,13$ need to be considered.
By applying Lemmata \ref{lem:autotype} to \ref{lem:none2},
all other automorphisms types $(p,f_U,f_V)$ except those
listed in Table \ref{tbl:types} are excluded. The last four
columns of the table---1F, V1, V2, and Seeds---are related 
to the main search to be discussed in the next section.

\begin{table}
\begin{center}
\begin{tabular}{ccccccr}\hline
\raisebox{-0.5mm}{$p$} & \raisebox{-0.5mm}{$f_U$} & \raisebox{-0.5mm}{$f_V$} &  \raisebox{-0.5mm}{1F} &  \raisebox{-0.5mm}{V1} &  \raisebox{-0.5mm}{V2} & \raisebox{-0.5mm}{Seeds} \\[0.5mm]\hline
\raisebox{-0.5mm}{2}   &   \raisebox{-0.5mm}{1}   &   \raisebox{-0.5mm}{2}   &  \raisebox{-0.5mm}{-}  &  \raisebox{-0.5mm}{F}  &  \raisebox{-0.5mm}{M}  & \raisebox{-0.5mm}{2579}  \\[0.5mm]
2   &   3   &   0   &  F  &  M  &  -  & 695   \\
2   &   3   &   4   &  -  &  F  &  M  & 10256 \\
2   &   5   &   0   &  F  &  M  &  -  & 894   \\
2   &   5   &   6   &  -  &  F  &  M  & 1206  \\
2   &   7   &   0   &  F  &  M  &  -  & 447   \\
3   &   1   &   2   &  F  &  F  &  F  & 65    \\
5   &   3   &   4   &  F  &  F  &  F  & 8     \\
7   &   6   &   0   &  M  &  M  &  M  & 9     \\
13  &   0   &   1   &  M  &  F  &  M  & 14    \\\hline
\end{tabular}
\end{center}
\bigskip
\caption{Automorphism types}
\label{tbl:types}
\end{table}

\subsection{The classification}

\label{sect:class}

We now present our approach for classifying the one-factorizations
that admit at least one automorphism of the types in Table~\ref{tbl:types}.
For this task we essentially rely on the framework established in 
\cite{K}.

The classification is based on certain substructures, to be
called \emph{seeds}, at least one of which is contained 
in every one-factorization that we want to classify. Due to the 
requirement of nontrivial symmetry, the precise definition of a
seed will unfortunately be somewhat technical. 
A simplified intuition to keep in mind is as follows. 
Consider an arbitrary one-factorization, $\mathcal{X}$, and 
(by some rule) select a set $T$ of points.
Let $\mathcal{S}$ be the set of blocks in $\mathcal{X}$ that have 
nonempty intersection with $T$. Now, based on the defining properties 
of one-factorizations (and the rule for selecting $T$), we 
can anticipate the structure of $\mathcal{S}$ without knowing
all the possible $\mathcal{X}$ explicitly, and classify 
all possible $\mathcal{S}$ up to isomorphism. 
These sets $\mathcal{S}$ will intuitively be the seeds;
the added technicality follows because 
(i) we insist on nontrivial symmetry in the form of a 
    prime-order group of automorphisms $\Pi\leq\Gamma$; and 
(ii) we must make precise the structure of $\mathcal{S}$ 
     in relation to $\Pi$ and $T$, without reference to any 
     containing one-factorization.

The technical definition of a seed is as follows.
Let $\Pi\leq\Gamma$ be a prime-order subgroup whose nonidentity 
elements have one of the types in Table~\ref{tbl:types}. 
Let $T\subseteq U\cup V$ be a set of points
such that the type of $\Pi$ and the columns 1F, V1, and V2 in 
Table~\ref{tbl:types} determine the size and composition of $T$ 
in relation to $\Pi$ as follows.
The size of $T$ is determined by the number of columns 
containing either an F (``fixed'') or an M (``moved'').
The column 1F (``one-factor'') indicates whether $T$
contains a point from $U$ and whether the point is fixed
or moved by $\Pi$.
The columns V1 (``first vertex'') and V2 (``second vertex'') indicate 
whether $T$ contains points from $V$ and whether these points are fixed
or moved by $\Pi$. (For example, for type $p=3$, $f_U=1$, $f_V=2$,
the set $T$ consists of one element of $U$ fixed by $\Pi$ and 
two elements of $V$ both moved by $\Pi$.)
Finally, let $\mathcal{S}$ be a union 
of $\Pi$-orbits of 3-subsets of $U\cup V$ such that, referring to
the elements of $\mathcal{S}$ as blocks,

\begin{itemize}
\item[(a')] for all $1\leq k\leq n-1$ and $1\leq i\leq n$,
the pair $\{u_k,v_i\}$ occurs in at most one block; 
\item[(b')] for all $1\leq i<j\leq n$, 
the pair $\{v_i,v_j\}$ occurs in at most one block; 
\item[(c')] for all $1\leq k<\ell\leq n-1$, no block contains the
pair $\{u_k,u_\ell\}$;
\item[(d')] for every $u_k\in T\cap U$, the point $u_k$ occurs in
exactly $n/2$ blocks;
\item[(e')] for every $v_i\in T\cap V$, the point $v_i$ occurs in
exactly $n-1$ blocks;
\item[(f')] the set $T$ has nonempty intersection with at least one
    block on every $\Pi$-orbit on $\mathcal{S}$; and
\item[(g')] the set $T$ occurs in at least one block.
\end{itemize}

Each tuple $(\Pi,T,\mathcal{S})$ meeting these requirements 
is called a \emph{seed}.
Two seeds, $(\Pi,T,\mathcal{S})$ and $(\Pi',T',\mathcal{S}')$,
are \emph{isomorphic} if there exists a $\gamma\in\Gamma$
such that $\gamma\Pi\gamma^{-1}=\Pi'$, $\gamma(T)=T'$, and
$\gamma(\mathcal{S})=\mathcal{S}'$. The permutation $\gamma$
is an \emph{isomorphism} of $(\Pi,T,\mathcal{S})$ onto 
$(\Pi',T',\mathcal{S}')$. An \emph{automorphism} of a seed is an isomorphism
of the seed onto itself. We denote by $\Aut(\Pi,T,\mathcal{S})$ the
automorphism group of a seed.
A one-factorization $\mathcal{X}$ \emph{contains} (or \emph{extends})
a seed $(\Pi,T,\mathcal{S})$ if $\Pi\leq\Aut(\mathcal{X})$ and
$\mathcal{S}\subseteq\mathcal{X}$.

We classify the seeds up to isomorphism by enlarging the
set $T$ one point at a time using the algorithms described
in \cite{K}. The number of nonisomorphic seeds associated with 
each automorphism type is displayed in the column Seeds in 
Table~\ref{tbl:types}.

Because every one-factorization with a nontrivial automorphism
group contains at least one seed, we can visit every isomorphism
class of one-factorizations by extending each classified seed
in all possible ways. The task of finding all one-factorizations
that contain a given seed $(\Pi,T,\mathcal{S})$ 
is an instance of the \emph{exact cover problem}.
Put otherwise, 
we must cover the remaining uncovered pairs of points of the form 
$\{u_k,v_i\}$ and $\{v_i,v_j\}$ exactly once in all possible ways 
using $\Pi$-orbits of triples of the form $\{u_k,v_i,v_j\}$;
each triple covers the pairs that occur in it. 
For algorithms, we refer to \cite{KP,Knuth00}.

We reject isomorphs among the visited one-factorizations 
using the framework in \cite{K}, which is an instantiation
of the canonical augmentation technique developed by McKay \cite{M3}. 
In essence, we identify a canonical $\Aut(\mathcal{X})$-orbit of
seeds contained by a visited one-factorization $\mathcal{X}$, and 
then check whether the seed $(\Pi,T,\mathcal{S})$ from which 
$\mathcal{X}$ was extended is in the canonical orbit. 
If yes, we accept $\mathcal{X}$ if it is also the (lexicographic) 
minimum of its $\Aut(\Pi,T,\mathcal{S})$-orbit.
Otherwise we reject $\mathcal{X}$.

In the search we find 581,042,656,543 one-factorizations that extend 
the classified seeds; among these we find 10,305,262,573 nonisomorphic 
one-factorizations with a nontrivial automorphism group.
Table \ref{tbl:autostat} displays the possible orders $i$ for 
the automorphism group and the associated number $N_i$ of 
nonisomorphic one-factorizations. ($N_1$ will be determined in 
Section \ref{sect:appl}.)
The search was distributed to a network of 180 Linux PCs
using the batch system \texttt{autoson} \cite{M2}, and
required in total a little over 5 years of CPU time.
The prevalent PC model in the network was Dell OptiPlex 745
with a 2.13-GHz Intel Core 2 Duo 6400 CPU and 2 GB main memory.

\begin{table}
\begin{center}
\begin{tabular}{rrrr}
\hline
\raisebox{-0.5mm}{$i$} & \raisebox{-0.5mm}{$N_i$} & \raisebox{-0.5mm}{$i$} & \raisebox{-0.5mm}{$N_i$} \\[0.5mm]\hline
\raisebox{-0.5mm}{1}   & \raisebox{-0.5mm}{\ \ 1132835411296799774} &\raisebox{-0.5mm}{21}  & \raisebox{-0.5mm}{1} \\[0.5mm]   
2   & 10300646080    	     &24     & 3 \\
3   & 4497762 	     	     &32     & 13 \\			  
4   & 104560 	     	     &39     & 3 \\			
5   & 2742 	     	     &42     & 2 \\			
6   & 9247 	     	     &48     & 1 \\	
8   & 1790 	     	     &64     & 3 \\
10  & 168 	     	     &84     & 1 \\
12  & 76 	     	     &156    & 1 \\
13  & 10 	     	     &192    & 1 \\			   
16  & 109 	     	     &{\ \ Total}  & 1132835421602062347  \\\hline
\end{tabular}
\end{center}
\bigskip
\caption{Nonisomorphic one-factorizations of $K_{14}$}
\label{tbl:autostat}
\end{table}

\subsection{Applying the Orbit-Stabilizer Theorem}

\label{sect:appl}

As a result of the computer searches, we know the value 
$\LF(K_{14})$ and the values $i$ and $N_i$ for each $i\geq 2$.
We now apply \eqref{eq:totaln} in the GDD representation
over fixed but arbitrary sets $U$ and $V$. 
Accordingly, $\Omega$ is the set of distinct one-factorizations
over $U$ and $V$, and $\Gamma$ is the group of all permutations
of $U\cup V$ that fix $U$ and $V$ setwise. Clearly, $|\Gamma|=13!\cdot 14!$\,.
Because there are $13!$ ways to label the one-factors in each 
distinct one-factorization of $K_{14}$ in the standard (graphical) 
representation, we have $|\Omega|=13!\,\cdot\,\LF(K_{14})$\,.
Solving \eqref{eq:totaln} for $N_1$ and computing $\sum_{i\geq 1} N_i$,
we have that the complete graph $K_{14}$ has exactly 
\[
\NF(K_{14})=\text{1,132,835,421,602,062,347} 
\]
nonisomorphic one-factorizations.

\subsection{Correctness checks}

\label{sect:correctness}

We carry out two checks to gain confidence in the correctness
of the values $N_i$ for $i\geq 2$ in Table \ref{tbl:autostat}.

First, both authors independently classified the seeds up to 
isomorphism, with identical results. One classification was 
conducted using the tools in \cite{K} and the other using a 
basic backtrack search with isomorph rejection based on recorded 
canonical forms. 

Second, for each automorphism type $(p,f_U,f_V)$, we count in two
different ways the distinct tuples $(\mathcal{X},\Pi,T,\mathcal{S})$, 
where $\mathcal{X}$ is a one-factorization and $(\Pi,T,\mathcal{S})$
is a seed of type $(p,f_U,f_V)$ contained in $\mathcal{X}$.

The double count is implemented as follows.
Fix an automorphism type $(p,f_U,f_V)$.
For a one-factorization $\mathcal{X}$, denote by
$\seeds(\mathcal{X})$ the number of distinct seeds of the fixed 
type that $\mathcal{X}$ contains. 
For a seed $(\Pi,T,\mathcal{S})$ of the fixed type, 
denote by $\ext(\Pi,T,\mathcal{S})$
the number of distinct one-factorizations that extend $(\Pi,T,\mathcal{S})$.
Taking the sum over all isomorphism classes of one-factorizations on 
the left-hand side, and over all isomorphism classes of seeds of 
the fixed type on the right-hand side, we have, by the Orbit-Stabilizer 
Theorem,
\begin{equation}
\label{eq:check}
\sum_{[\mathcal{X}]}
\frac{|\Gamma|}{|\Aut(\mathcal{X})|}\cdot\seeds(\mathcal{X})=
\sum_{[(\Pi,T,\mathcal{S})]}
\frac{|\Gamma|}{|\Aut(\Pi,T,\mathcal{S})|}\cdot\ext(\Pi,T,\mathcal{S}).
\end{equation}

To evaluate the right-hand side of \eqref{eq:check}, we 
record the number of extensions $\ext(\Pi,T,\mathcal{S})$ 
and $|\Aut(\Pi,T,\mathcal{S})|$ for each classified seed $(\Pi,T,\mathcal{S})$.
In particular, $\ext(\Pi,T,\mathcal{S})$ is simply the number of solutions
found in the search for exact covers. The right-hand sides 
obtained for each automorphism type are listed in column Count 
in Table \ref{tbl:check}.

\begin{table}
\begin{center}
\begin{tabular}{cccrr}\hline
\raisebox{-0.5mm}{$p$} & \raisebox{-0.5mm}{$f_U$} & \raisebox{-0.5mm}{$f_V$} & \raisebox{-0.5mm}{$m(p,f_U,f_V)$} & \raisebox{-0.5mm}{Count}\\[0.5mm]\hline
\raisebox{-0.5mm}{2}   &   \raisebox{-0.5mm}{1}   &   \raisebox{-0.5mm}{2}   & \raisebox{-0.5mm}{24} & \raisebox{-0.5mm}{598566905953570569439936512000}  \\[0.5mm]
2   &   3   &   0   &  42  & 10362562621908790673701601280000      \\
2   &   3   &   4   &  40  & 11064634216108608459689164800000      \\
2   &   5   &   0   &  70  & 109764651070947200382428774400000     \\
2   &   5   &   6   &  48  & 314439007330643189170176000000        \\
2   &   7   &   0   &  98  & 65958792604977492862823301120000      \\
3   &   1   &   2   &  1   & 814728009186504568995840000           \\
5   &   3   &   4   &  6   & 1840950333789938122752000             \\
7   &   6   &   0   &  49  & 2850020421206999040000                \\
13  &   0   &   1   &  13  & 6016709778103664640000                \\\hline
\end{tabular}
\end{center}
\bigskip
\caption{Double counting check}
\label{tbl:check}
\end{table}

The left-hand side of \eqref{eq:check} is accumulated
for each classified one-factorization, $\mathcal{X}$.
For each such $\mathcal{X}$, we find all distinct prime-order subgroups
$\Pi\leq\Aut(\mathcal{X})$. For each such $\Pi$, we find its type
$(p,f_U,f_V)$, and accumulate the left-hand side of \eqref{eq:check}
for this type by $|\Gamma|/|\Aut(\mathcal{X})|\cdot m(p,f_U,f_V)$, 
where $m(p,f_U,f_V)$ is the number of distinct seeds 
of the form $(\Pi,T,\mathcal{S})$ contained in $\mathcal{X}$. 
Because $\mathcal{S}$ is uniquely determined by $\mathcal{X}$, $\Pi$, 
and $T$, we can determine $m(p,f_U,f_V)$ by combinatorial arguments based 
on $(p,f_U,f_V)$ and Table~\ref{tbl:types}.
For example, consider $p=5$, $f_U=3$, $f_V=4$.
From Table \ref{tbl:types} we find that $T$ consists
of one point of $U$ fixed by $\Pi$ and two points
of $V$ fixed by $\Pi$. By (g') and $\mathcal{S}\subseteq\mathcal{X}$,
we have that $T$ is a block of $\mathcal{X}$.
There are $\binom{4}{2}=6$ possibilities to select a pair of points 
in $V$ fixed by $\Pi$, each of which occurs in a unique block
of $\mathcal{X}$. Thus, $m(5,3,4)=6$. The other values
$m(p,f_U,f_V)$ are displayed in Table \ref{tbl:check}.

For each type $(p,f_U,f_V)$, we find that the computed left-hand 
and right-hand sides of \eqref{eq:check} agree. This, together
with the observation that the left-hand side of \eqref{eq:check}
depends on each of the computed values $i$ and $N_i$ for $i\geq 2$, 
gives us confidence that Table~\ref{tbl:autostat} is correct.

\section*{Acknowledgments}
This research was supported in part by 
the Academy of Finland, Grants 107493, 110196, and 117499.

%%%%%%%%%%%%%%%%%%%%%%%%%%%%%%%%%%%%%%%%%%%%%%%%%%%%%%%%%%%%%%%% References %%%


\begin{thebibliography}{99}
\bibitem{A} L.~D.~Andersen, ``Factorizations of graphs,''
  Handbook of combinatorial designs, C.~J.~Colbourn and
  J.~H.~Dinitz (Editors), 2nd ed., 
  Chapman \& Hall/CRC, Boca Raton, 2007, pp.~740--755.
\bibitem{DS} L.~E.~Dickson and F.~H.~Safford, Solution to problem 8 
  (group theory), Amer Math Monthly 13 (1906), 150--151.
\bibitem{DGM} J.~H.~Dinitz, D.~K.~Garnick, and B.~D.~McKay,
  There are $526{,}915{,}620$ nonisomorphic one-factorizations
  of $K_{12}$, J Combin Des 2 (1994), 273--285.
\bibitem{GO}
  E.~N.~Gelling and R.~E.~Odeh,
  On $1$-factorizations of the complete graph and the relationship
  to round-robin schedules, Congr Numer 9 (1974), 213--221.
\bibitem{gmplib}
  GNU Multiple Precision Arithmetic Library, Version 4.2.1.
  Available at $\langle$\url{http://gmplib.org/}$\rangle$.
\bibitem{IP} E.~Ihrig and E.~Petrie,
  Automorphism groups of $1$-factorizations of $K_{12}$,
  Congr Numer 109 (1995), 179--192.
\bibitem{K} P.~Kaski, Isomorph-free exhaustive generation of designs 
  with prescribed groups of automorphisms, SIAM J Discrete Math
  19 (2005), 664--690.
\bibitem{KO} P.~Kaski and P.~R.~J.~{\"O}sterg{\aa}rd, 
  Classification algorithms for codes and designs,
  Springer, Berlin, 2006.
\bibitem{KP} P.~Kaski and O.~Pottonen, 
  {\sf libexact} User's Guide, Version 1.0, in preparation.
\bibitem{Knuth00} D.~E.~Knuth, ``Dancing links,''
  in Millennial perspectives in computer science,
  J.~Davies, B.~Roscoe, and J.~Woodcock (Editors),
  Palgrave Macmillan, Basingstoke, 2000, pp.~187--214.
\bibitem{M1} B.~D.~McKay, \emph{nauty} user's guide 
  (version\/ $1.5$), Technical Report TR-CS-90-02, 
  Computer Science Department, Australian National University, 
  Canberra, 1990.
\bibitem{M2} B.~D.~McKay,
  \texttt{autoson} -- a distributed batch system for UNIX
  workstation networks (version 1.3),
  Technical Report TR-CS-96-03,
  Computer Science Department, Australian National University, 1996.
\bibitem{M3} B.~D.~McKay, 
  Isomorph-free exhaustive generation,
  J Algorithms 26 (1998), 306--324.
\bibitem{MMM} B.~D.~McKay, A.~Meynert, and W.~Myrvold, 
  Small Latin squares, quasigroups,
  and loops, J Combin Des 15 (2007), 98--119.
\bibitem{MW} B.~D.~McKay and I.~M.~Wanless, On the number of Latin 
  squares, Ann Comb 9 (2005) 335-344. 
\bibitem{Me} M.~Meringer, Fast generation of regular graphs and 
  construction of cages, J. Graph Theory 30 (1999), 137--146.
\bibitem{SS} E. Seah and D. R. Stinson,
  On the enumeration of one-factorizations of complete
  graphs containing prescribed automorphism groups,
  Math Comp 50 (1988), 607--618.
\bibitem{W1} W.~D.~Wallis,
  ``One-factorizations of complete graphs,''
  in Contemporary design theory: A collection of surveys,
  J.~H.~Dinitz and D.~R.~Stinson (Editors),
  Wiley, New York, 1992, pp.~593--631.
\bibitem{W2} W.~D.~Wallis,
  One-factorizations, Kluwer, Dordrecht, 1997.
\bibitem{W3}  D.~B.~West, Introduction to graph theory,
  2nd ed., Prentice--Hall, Upper Saddle River, 2001.

\end{thebibliography}
\end{document}